\documentclass[12pt,
]{amsart}

\usepackage{amssymb,amsfonts,amstext,amsmath}
\usepackage{amsthm,latexsym,mathrsfs}
\usepackage{mathbbol}

\usepackage{mparhack}
\usepackage{color}

\theoremstyle{definition}

\theoremstyle{remark}

\numberwithin{equation}{section}

\newcommand{\abs}[1]{\lvert#1\rvert}

\input{gklm.sty}

\begin{document}

\title{A Groszek - Laver pair of undistinguishable
 $\Eo$ classes}

\author{Mohammad~Golshani}
\address{School of Mathematics, Institute for Research in
Fundamental Sciences (IPM), P.O. Box:
19395-5746, Tehran-Iran}
\email{golshani.m@gmail.com}
\thanks{The first author was supported in part by
IPM Grant \#91030417.}
\author{Vladimir~Kanovei}
\address{IITP RAS and MIIT,
  Moscow, Russia}
\email{kanovei@googlemail.com}
\thanks{The second author was supported in part by
RFBR Grant \#13-01-00006.}
\author{Vassily~Lyubetsky}
\address{IITP RAS,
  Moscow, Russia}
\email{lyubetsk@iitp.ru}
\thanks{The second and third authors were supported in part by
RNF Grant \#14-50-00150.}

\subjclass[2000]{Primary 03E15, 03E35; Secondary 03E45}

\date{March 07, 2015}

\keywords{Forcing, equivalence classes, ordinal definability,
Groszek - Laver pair}

\begin{abstract}
A generic extension $\rL[x,y]$ of $\rL$ by reals $x,y$ is defined,
in which the union of $\mathsf E_0$-classes of $x$ and $y$
is a $\ip12$ set, but neither of these two $\mathsf E_0$-classes
is separately ordinal-definable.
\end{abstract}

\maketitle

\section{Introduction}
\las{cha1}

Let a \rit{Groszek - Laver pair} be any unordered $\od$
(ordinal-definable) pair
$\ans{X,Y}$ of sets $X,Y\sq\bn$ such that neither of $X,Y$ is
separately $\od$.
As demonstrated in \cite{gl}, if $\ang{x,y}$ is a Sacks$\ti$Sacks
generic pair of reals over $\rL$, the constructible universe,
then their degrees of
constructibility $X=[x]_\rL\cap\bn$ and $Y=[y]_\rL\cap\bn$
form such a pair in $\rL[x,y]$;
the set $\ans{X,Y}$ is definable as
the set of all \dd\rL degrees of reals,
\dd\rL minimal over $\rL$.

As the sets $X,Y$ in this example are obviously uncountable,
one may ask whether there can consistently exist
a Groszek -- Laver pair of \rit{countable} sets.
The next theorem answers this question in the positive in
a rather strong way: both sets are \dd\Eo classes in the
example!
(Recall that the \eqr\ $\Eo$ is defined on $\dn$ as follows:
$x\Eo y$ iff $x(n)=y(n)$ for all but finite $n$.)

\bte
\lam{mt}
It is true in a suitable generic extension\/ $\rL[x,y]$
of\/ $\rL$, by a pair of reals\/
$x,y\in\dn$ that the union of\/ \dd\Eo equivalence classes\/
$\eko x\cup\eko y$
is\/ $\ip12$, but neither of the sets\/ $\eko x,\eko y$
is separately $\od$.
\ete

The forcing we employ is a conditional product $\dPd$ of an
``\dd\Eo large tree''\snos
{An \dd\Eo large tree is a perfect tree $T\sq\bse$ such that
${\Eo}\res{[T]}$ is not smooth, see \cite[10.9]{kanB}.}
version $\dP$  of a forcing notion,
introduced in \cite{kl:ceo} to define a model
with a $\ip12$ \dd\Eo class containing no $\od$ elements.
The forcing in \cite{kl:ceo} was a clone of Jensen's minimal
$\ip12$ real singleton forcing \cite{jenmin}
(see also Section 28A of \cite{jechmill}),
but defined on the base of the
Silver forcing instead of the Sacks forcing.
The crucial advantage of Silver's forcing here is that it
leads to a Jensen-type forcing naturally closed under the 0-1
flip at any digit, so that the corresponding extension contains
a $\ip12$ \dd\Eo class of generic reals instead of a
$\ip12$ generic singleton as in \cite{jenmin}.

In another relevant note \cite{kl:cds} it is demonstrated
that a countable $\od$ set of reals (not an \dd\Eo class),
containing no $\od$ elements, exists in a generic extension
of $\rL$ via the countable finite-support product of
Jensen's \cite{jenmin} forcing itself.
The existence of such a set was discussed as an open
question at the \rit{Mathoverflow} website\snos
{\label{snos1}
A question about ordinal definable real numbers.
\rit{Mathoverflow}, March 09, 2010.
{\tt http://mathoverflow.net/questions/17608}.
}
and at FOM\snos
{\label{snos2}%
Ali Enayat. Ordinal definable numbers. FOM Jul 23, 2010.
{\tt http://cs.nyu.edu/pipermail/fom/2010-July/014944.html}}%
,
and the result in \cite{kl:cds} was
conjectured by Enayat (Footnote~\ref{snos2})
on the base of his study of finite-support products of
Jensen's forcing in \cite{ena}.

The remainder of the paper is organized as follows.

We introduce \dd\Eo large perfect trees in $\bse$ in
Section~\ref{tre},
study their splitting properties in Section~\ref{sst},
and consider \dd\Eo large-tree forcing notions
in Section~\ref{sitf}, \ie, collections of
\dd\Eo large trees closed under both restriction and
action of a group of transformations naturally associated
with $\Eo$.

If $\dP$ is an \dd\Eo large-tree forcing notion then the
\rit{conditional product forcing} $\dpd$ is a part of the
full forcing product $\dP\ti\dP$ which
contains all conditions $\ang{T,T'}$ of trees $T,T'\in\dP$,
\dd\Eo connected in some way.
This key notion, defined in Section~\ref{cop}, goes back
to early research on the Gandy -- Harrington forcing
\cite{hms,hkl}.

The basic \dd\Eo large-tree forcing $\dP$ employed in the
proof of Theorem~\ref{mt} is defined, in $\rL$, in the form
$\dP=\bigcup_{\xi<\omi}\dU_\xi$ in Section~\ref{jfor}.
The model $\rL[x,y]$ which proves the theorem is then a
\dd{(\dpd)}generic extension of $\rL$; it is studied in
Section~\ref{mod}.
The elements $\dU_\xi$ of this inductive construction are
countable \dd\Eo large-tree forcing notions in $\rL$.

The key issue is, given a subsequence
$\sis{\dU_\eta}{\eta<\xi}$ and accordingly the union
$\dP_{<\xi}=\bigcup_{\eta<\xi}{\dU_\eta}$, to define the
next level $\dU_\xi$.
We maintain this task in Section~\ref{jex} with the help of
a well-known splitting/fusion construction, modified so
that it yields \dd\Eo large perfect trees.
Generic aspects of this construction lead to the CCC property
of $\dP$ and $\dpd$ and very simple reading of real names,
but most of all to the crucial property that if $\ang{x,y}$
is a pair of reals \dd{(\dpd)}generic over $\rL$ then any
real $z\in\rL[x,y]$ \dd\dP generic over $\rL$ belongs to
$\eko x\cup \eko y$.
This is Lemma~\ref{only} proved, on the base of preliminary
results in Section~\ref{saway2}.

The final Section~\ref{konk} briefly discusses some related
topics.

\section{\dd\Eo large trees}
\las{tre}

Let $\bse$ be the set of all strings (finite sequences)
\index{string}%
\index{string!empty, $\La$}%
\index{zz2.omega@$\bse$}%
of numbers $0,1$,
including the empty string $\La$.
If $t\in\bse$ and $i=0,1$ then
$t\we i$ is the extension of $t$ by $i$ as the rightmost term.
\index{extension!t.k@$t\we k$}%
\index{extension!s.t@$s\we t$}%
\index{zzt.k@$t\we k$}%
\index{zzs.t@$s\we t$}%
If $s,t\in\bse$ then $s\sq t$ means that $t$ extends $s$,
\index{zzssubseteqt@$s\sq t$}%
\index{zzssubsett@$s\su t$}%
$s\su t$ means proper extension, and $s\we t$ is the
concatenation.
If $s\in\bse$ then $\lh s$ is the length of $s$,
\index{length, $\lh s$}%
and we let
$2^n=\ens{s\in\bse}{\lh s=n}$ (strings of length $n$).%
\index{string!of length $n$}%

Let any $s\in\bse$ {\ubf act} on $\dn$ so that
\index{action!s*t@$s\app t$}%
\index{action!s*X@$s\app X$}%
\index{zzs*t@$s\app t$}%
\index{zzs*X@$s\app X$}%
$(s\app x)(k)=x(k)+s(k)\pmod 2$ whenever $k<\lh s$ and simply
$(s\app x)(k)=x(k)$ otherwise.
If $X\sq\dn$ and $s\in\bse$ then, as usual, let
$s\app X=\ens{s\app x}{x\in X}$.

Similarly if $s,t\in \bse$ and $\lh s=m\le n=\lh t$, then define
$s\app t\in 2^n$ so that
\index{zz*@$\app\;$}%
$(s\app t)(k)=t(k)+s(k)\pmod 2$ whenever $k<m$ and
$(s\app t)(k)=t(k)$ whenever $m\le k<n$.
If $m>n$ then let simply $s\app t=(s\res n)\app t$.
Note that $\lh{s\app t}=\lh t$ in both cases.
Let $s\app T=\ens{s\app t}{t\in T}$ for $T\sq\bse$.

If $T\sq\bse$ is a tree and $s\in T$ then put
$\req Ts=\ens{t\in T}{s\sq t\lor t\sq s}$.
\index{restriction!T:s@$\req Ts$}%
\index{zzT:s@$\req Ts$}%

Let $\pet$ be the set of all \rit{perfect} trees
\index{tree!perfect}%
\index{zzPT@$\pet$}%
$\pu\ne T\sq \bse$
\imar{pet}%
(those with no endpoints and no isolated branches).
If $T\in\pet$ then there is a largest string $s\in T$
such that $T=T\ret s$; it is denoted by $s=\roo T$
\index{stem, $\roo T$}%
\index{zzstemT@$\roo T$}%
(the {\it stem\/} of $T$);
we have $s\we 1\in T$ and $s\we 0\in T$ in this case.
If $T\in\pet$ then
$$
[T]=\ens{a\in\dn}{\kaz n\,(a\res n\in T)}\sq\dn
$$
\index{zz:T:@$[T]$}%
is the perfect set of all \rit{paths through $T$};
clearly $[S]\sq [T]$ iff $S\sq T$.

Let $\pes$ (large trees)
\index{tree!special large}%
\index{zzLT@$\pes$}%
be the set of all \rit{\elt s}: those
$T\in\pet$ such that there is a
double sequence of non-empty
strings $\ee in=\ee in(T)\in\bse$,
\index{zzqinT@$\ee in(T)$}%
$n<\om$ and $i=0,1$, such that\vom
\bit
\item[$-$]
$\lh{\ee0n}=\lh{\ee1n}\ge1$ and $\ee in(0)=i$ for all $n$;
\vom

\item[$-$]
$T$ consists of all substrings
of strings of the form
$r\we \ee{i(0)}0\we \ee{i(1)}1\we \dots\we \ee{i(n)}n$
in $\bse,$
where $r=\roo T$, $n<\om$, and $i(0),i(1),\dots,i(n)\in\ans{0,1}$.
\vom
\eit

\noi
We let $\oin T0=\lh r$ and then by induction
$\oin T{n+1}=\oin Tn+\lh{\ee in}$, so that
\index{zzsplnT@$\oin Tn$}%
$\oi T=\ens{\oin Tn}{n<\om}\sq\om$ is the set of
\index{zzsplT@$\oi T$}%
\rit{splitting levels} of $T$.
Then
$$
[T]=\ens{a\in\dn}
{a\res \lh r=r\land
\kaz n\,\big(a\res[\oin Tn,\oin T{n+1})= \ee0n\,\text{ or }\,\ee1n)}.
$$

\ble
\lam{trans}
Assume that\/ $T\in\pes$ and\/ $h\in\oi T$.
Then\/\vom
\ben
\renu
\itla{trans1}
if\/ $u,v\in2^h\cap T$ then\/
$\req Tv=(u\app v)\app\req Tu$ and\/ $(u\app v)\app T=T\;;$\vom

\itla{trans2}
if\/ $\sg\in\bse$ then\/
$T=\sg\app T$ or\/ $T\cap(\sg\app T)$ is finite.\vom
\een
\ele
\bpf
\ref{trans2}
Suppose that $T\cap(\sg\app T)$ is infinite.
Then there is an infinite branch $x\in[T]$ such that
$y=\sg\app x\in [T]$, too.
We can assume that $\lh\sg$ is equal to some $h=\oin Tn$.
(If $\oin T{n-1}<h<\oin Tn$ then extend $\sg$ by
$\oin Tn-h$ zeros.)
Then $\sg=(x\res h)\app(y\res h)$.
It remains to apply \ref{trans1}.
\epf

\bpri
\lam{clop}
If $s\in\bse$ then $T[s]=\ens{t\in\bse}{s\sq t\lor t\su s}$
\index{zzT:s:@$T[s]$}%
is a tree in $\pes$, $\roo{T[s]}=s$, and\/ $\ee in(T[s])=\ang i$
for all $n,i$.
Note that $T[\La]=\bse$ (the full binary tree), and
$\req{T[\La]}s=\req{(\bse)}s=T[s]$ for all $s\in\bse$.
\epri

\section{Splitting of large trees}
\las{sst}

The \rit{simple splitting\/} of a tree $T\in\pes$ consists
\index{splitting}%
of smaller trees
$$
\raw T0=T\ret{\roo T\we 0}
\quad\text{and}\quad
\raw T1=T\ret{\roo T\we 1}\,,
$$
\index{zzT->i@$\raw T i$}%
so that $[\raw Ti]=\ens{x\in[T]}{x(h)=i}$,
where $h=\oin T0=\lh{\roo T}$.
Clearly $\raw Ti\in\pes$ and
$\oi{\raw Ti}= \oi T\bez\ans{\oin T0}$.

\ble
\lam{uv}
If\/ $R,S,T\in\pes$, $S\sq \raw R0$, $T\sq \raw R1$,
$\sg\in\bse,$ $T=\sg\app S$, and\/
$\lh\sg\le \lh{\roo S}=\lh{\roo T}$
then\/ $U=S\cup T\in\pes$, $\roo U=\roo R$, and\/
$S=\raw U0$, $T=\raw U1$.
\qed
\ele

The splitting can be iterated, so that if
$s\in 2^n$ then we define
$$
\raw Ts=\raw{\raw{\raw{\raw T{s(0)}}{s(1)}}{s(2)}\dots}{s(n-1)}
\,.
$$
\index{zzT->s@$\raw T s$}%
We separately define $\raw T\La=T$, where $\La$ is the empty
string as usual.

\ble
\lam{cloL}
In terms of Example~\ref{clop}, $T[s]=\raw{(\bse)}s=\req{(\bse)}s$,
$\kaz s$.
Generally if\/ $T\in\pes$ and\/ $2^n\sq T$ then\/
$\raw Ts=\req Ts$ for all\/ $s\in 2^n$.\qed
\ele

If $T,S\in\pes$ and $n\in\om$ then let $S\nq n T$
($S$ \rit{\dd nrefines\/ $T$})
\index{zzSsubnT@$S\nq n T$}%
\index{zzsubn@$\nq n$}%
mean that
$S\sq T$ and $\oin Tk=\oin Sk$ for all $k<n$.
In particular, $S\nq 0 T$ iff simply $S\sq T$.
By definition if $S\nq{n+1} T$ then $S\nq n T$
(and $S\sq T$), too.

\ble
\lam{sadd}
Suppose that\/ $T\in\pes$, $n<\om$, and\/ $h=\oin Tn$.
Then\vom
\ben
\renu
\itla{sadd1}
$T=\bigcup_{s\in2^n}\raw Ts$ and\/
\imar{sadd1}
$[\raw Ts] \cap[\raw Tt]=\pu$ for all\/ $s\ne t$ in\/ $2^n\;;$\vom


\itla{sadd2}
if\/ $S\in\pes$ then\/ $S\nq nT$
\imar{sadd2}
{\ubf iff\/} $\raw Ss\sq\raw Ts$ for all strings\/
$s\in 2^{\leqs n}$
{\ubf iff\/} $S\sq T$ and\/ $S\cap 2^h=T\cap 2^h\;;$\vom

\itla{sadd3}
if\/ $s\in2^n$ then\/ $\lh{\roo{\raw Ts}}=h$
and there is a string\/
\imar{sadd3}
$u[s]\in 2^h\cap T$ such that\/
$\raw Ts=\req T{u[s]}\;;$\vom

\itla{sadd4}
if\/ $u\in 2^h\cap T$
\imar{sadd4}
then there is a string\/ $s[u]\in2^n$ \st\/
$\req Tu=\raw T{s[u]}\;;$\vom

\itla{sadd5}
if\/ $s_0\in 2^n$ and\/ $S\in\pes$, $S\sq \raw T{s_0}$, then
\imar{sadd5}
there is a unique tree\/ $T'\in\pes$ such that\/
$T'\nq n T$ and\/ $\raw{T'}{s_0}=S\;.$
\een
\ele

\bpf
\ref{sadd3}
Define
$u[s]=\roo T\we
\ee{s(0)}0(T)\we \ee{s(1)}1(T)\we \dots\we \ee{s(n-1)}{n-1}(T)$.\vom
\index{zzu.s@$u[s]$}%
\iman

\ref{sadd4}
Define $s=s[u]\in2^n$ by $s(k)=u(\oin Tk)$ for all $k<n$.\vom
\index{zzs.u@$s[u]$}%

\ref{sadd5}
Let $u_0=u[s_0]\in 2^h$.
Following Lemma~\ref{trans},
define $T'$ so that $T'\cap2^h=T\cap2^h$, and if
$u\in T\cap2^h$ then $\req {T'}u=(u\cdot u_0)\cdot S$;
in particular $\req{T'}{u_0}=S$.
%
\epf

\ble
[fusion]
\lam{fus}
Suppose that\/
$\dots \nq 5 T_4\nq 4 T_3\nq 3 T_2\nq 2 T_1\nq 1 T_0$ is
an infinite decreasing sequence of trees in\/ $\pes$.
Then\/\vom
\ben
\renu
\itla{fus1}
$T=\bigcap_nT_n\in\pes\;;$\vom
\itla{fus2}
 if\/ $n<\om$ and\/ $s\in2^{n+1}$ then\/
$\raw Ts=T\cap\raw{T_{n}}s=\bigcap_{m\ge n}\raw{T_m}s$.
\een
\ele
\bpf
Both parts are clear, just note that  $\oi T=\ens{\oin{T_n}{n}}{n<\om}$.
\epf

\section{Large-tree forcing notions}
\las{sitf}

Let a {\it large-tree forcing notion\/}
\index{large-tree forcing notion}%
\index{forcing!large-tree forcing}%
\index{LTF}%
(\sif)
be any set
$\dP\sq\pes$ such that\vom
\ben
\snenu
\itla{sitf2}
\label{utp}
if $u\in T\in\dP$ then $T\ret u\in \dP$;\vom
\imar{sitf2}

\itla{sitf3}
if $T\in\dP$ and $s\in\bse$ then $s\app T\in \dP$.\vom
\imar{sitf3}
%
\een
We'll typically consider \sif s $\dP$
containing the full tree $\bse$.
In this case, $\dP$ contains all trees $T[s]$ of
Example~\ref{clop} by Lemma~\ref{cloL}.

Any \sif\ $\dP$ can be viewed as a forcing notion
(if $T\sq T'$ then $T$ is a stronger condition), and then
it adds a real in $\dn$.

If $\dP\sq\pes$, $T\in\pes$, $n<\om$, and all split
trees $\raw Ts$, $s\in 2^n$, belong to $\dP$, then we say
that $T$ is an \rit{\dd ncollage over\/ $\dP$}.
\index{collage@\dd ncollage}%
Let $\sct\dP n$ be the set of all trees $T\in\pes$ which are
\index{zzLCnP@$\sct\dP n$}%
\imar{sct Pn}
\dd ncollages over\/ $\dP$, and
$\sco\dP=\bigcup_n\sct \dP n$.
\index{zzLCP@$\sco\dP$}%
Note that $\sct\dP n\sq \sct\dP{n+1}$ by \ref{sitf2}.

\ble
\lam{stf}
Assume that\/ $\dP\sq\pes$ is a\/ \sif\ and $n<\om$.
Then\vom
\ben
\renu
\itla{stf1}
if\/ $T\in\pes$ and\/ $s_0\in2^n$
then\/ $\raw T{s_0}\in\dP$ iff\/ $T\in\sct\dP n\;;$\vom

\itla{stf2}
if\/ $P\in\sct\dP n$, $s_0\in2^n$, $S\in\dP$, and\/
$S\sq \raw P{s_0}$, then there is a tree\/ $Q\in\sct\dP n$
such that\/ $Q\nq n P$ and\/ $\raw{Q}{s_0}=S\;;$\vom

\itla{stf2c}
if\/ $P\in\sct\dP n$ and a set\/ $D\sq\dP$ is open dense
\imar{stf2c}
\iman
in\/ $\dP$, then there is a tree\/ $Q\in\sct\dP n$
such that\/ $Q\nq n P$ and\/ $\raw{Q}{s}\in D$
for all\/ $s\in 2^n\;;$\vom

\itla{stf3}
if\/ $P\in\sct\dP n$,
\imar{stf3}
$S,T\in\dP$, $s,t\in2^n,$ $S\sq\raw P{s\we0}$, $T\sq\raw P{t\we1}$,
$\sg\in\bse$, and\/ $T=\sg\app S$, then
there is a tree\/ $Q\in\sct\dP{n+1}$, $Q\nq{n+1}P$,
such that\/ $\raw{Q}{s\we0}\sq S$ and\/ $\raw{Q}{t\we1}\sq T$.
\een
\ele

Recall that a set\/ $D\sq\dP$ is \rit{open dense} in $\dP$ iff,
1st, if $S\in\dP$ then there is a tree $T\in D$, $T\sq S$,
and 2nd, if $S\in\dP$, $T\in D$, and $S\sq T$, then $S\in D$, too.

\bpf
\ref{stf1}
If $T\in\sct\dP n$ then by definition $\raw T{s_0}\in\dP$.
To prove the converse,
let $h=\oin Tn$, and let $h[s]\in 2^h\cap T$ satisfy
$\raw Ts=\req T{u[s]}$ for all $s\in2^n$ by
Lemma~\ref{sadd}\ref{sadd3}.
If $\raw T{s_0}\in\dP$ then
$\raw Ts=\req T{u[s]}=(u[s]\app u[s_0])\app \req T{u[s]}$
by Lemma~\ref{trans}, so $\raw Ts\in\dP$ by \ref{sitf3}.
Thus $T\in \sct\dP n$.\vom

\ref{stf2}
By Lemma~\ref{sadd}\ref{sadd5}
there is a tree $Q\in\pes$
such that $Q\nq n P$ and $\raw{Q}{s_0}=S$.
We observe that $Q$ belongs to $\sct\dP n$ by \ref{stf1}.\vom

\ref{stf2c}
Apply \ref{stf2} consecutively $2^n$ times (all $s\in2^n$).\vom

\ref{stf3}
We first consider the case when $t=s$.
If $\lh\sg\le L=\lh{\roo{S}}=\lh{\roo T}$ then by Lemma~\ref{uv}
$U=S\cup T\in\pes$, $\roo U=\roo{\raw Ps}$,
and $\raw U0=S$, $\raw U1=T$.
%
Lemma~\ref{sadd}\ref{sadd5} yields a tree $Q\in\pes$ such
that $Q\nq n P$ and $\raw Qs=U$, hence
$\roo{\raw Qs}=\roo{\raw Ps}$ by the above.
This implies $\oin Qn=\oin Pn$ by Lemma~\ref{sadd}\ref{sadd3},
and hence $Q\nq {n+1} P$.
\vyk{
$\raw Q{s\we0}=S\sq \raw P{s\we0}$,
$\raw Q{s\we1}=T\sq \raw P{s\we1}$.

We claim that generally
$\raw Q{s'\we i}\sq \raw P{s'\we i}$ for all $s'\in2^n$ and
$i=0,1$.
Indeed let $h=\oin Pn$.
Let $u=u[s]$ and $u'=u[s']$ be the strings in $P\cap 2^h$
defined by Lemma~\ref{sadd}\ref{sadd3}.

$\raw Q{s'\we i}=$

We conclude that $Q\nq{n+1} P$.
}%
And finally $Q\in \sct\dP{n+1}$ by \ref{stf1} since
$\raw Q{s\we0}=S\in\dP$.

Now suppose that $\lh\sg>L$.
Take any string $u\in S$ with $\lh u\ge\lh s$.
The set $S'=\req Su\sq S$ belongs to $\dP$ and obviously
$\lh{\roo{S'}}\ge \lh\sg$.
It remains to follow the case already considered for the trees
$S'$ and $T'=\sg\app S'$.

Finally consider the general case $s\ne t$.
Let $h=\oin Pn$, $H=\oin P{n+1}$.
Let $u=u[s]$ and $v=u[t]$ be the strings in $P\cap 2^h$
defined by Lemma~\ref{sadd}\ref{sadd3} for $P$, so that
$\req Pu=\raw Ps$ and $\req Pv=\raw Pt$, and let
$U,V\in2^H\cap P$ be defined accordingly so that
$\req PU=\raw P{s\we1}$ and $\req PV=\raw P{t\we1}$.
Let $\rho=u\app v$.
Then $\raw P{s}=\rho\app\raw P{t}$ by Lemma~\ref{trans}.
However we have $U=u\we\tau$ and $V=v\we\tau$ for one
and the same string $\tau$, see the proof of
Lemma~\ref{sadd}\ref{sadd3}.
Therefore $U\app V=u\app v=\rho$ and
$\raw P{s\we1}=\rho\app\raw P{t\we1}$
still by Lemma~\ref{trans}.

It follows that the tree $T_1=\rho\app T$ satisfies
$T_1\sq\raw P{s\we1}$.
Applying the result for $s=t$, we get a tree
$Q\in\sct\dP{n+1}$, $Q\nq{n+1}P$, such that\/
$\raw{Q}{s\we0}\sq S$ and\/ $\raw{Q}{s\we1}\sq T_1$.
Then by definition $\oin Pk=\oin Qk$ for all $k\le n$,
and $\raw Qs\sq\raw Ps$ for all $s\in2^{n+1}$ by
Lemma~\ref{sadd}\ref{sadd2}.
Therefore the same strings $u,v$ satisfy
$\req Qu=\raw Q{s}$ and $\req Qv=\raw Q{t}$.
The same argument as above implies
$\raw{Q}{t\we1}=\rho\app \raw{Q}{s\we1}$.
We conclude that
$\raw{Q}{t\we1}\sq \rho\app T_1=T$, as required.
\epf

\section{Conditional product forcing}
\las{cop}

Along with any \sif\ $\dP$, we'll consider the
{\ubf conditional product} $\dpd$, which by definition
\index{zzP2@$\dpd$}%
\index{forcing!P2@$\dpd$}%
\index{forcing!conditional product}%
consists of all pairs $\ang{T,T'}$ of trees $T,T'\in\dP$ such that
there is a string $s\in\bse$ satisfying $s\app T=T'$.
We order $\dpd$
componentwise so that $\ang{S,S'}\le\ang{T,T'}$
($\ang{S,S'}$ is stronger) iff $S\sq T$ and $S'\sq T'$.\snos
{Conditional product forcing notions of this kind were considered
in \cite{hms,hkl,k:sol} and some other papers
with respect to the Gandy -- Harrington and similar forcings,
and recently in
\cite{ksz} with respect to many forcing notions.}

\bre
\lam{ref1}
$\dpd$ forces a pair of \dd\dP generic reals.
Indeed if $\ang{T,T'}\in\dpd$ with $s\app T=T'$ and $S\in\dP$,
$S\sq T$, then there is a tree $S'=s\app S\in \dP$
(we make use of \ref{sitf3})
such that $\ang{S,S'}\in\dpd$
and  $\ang{S,S'}\le\ang{T,T'}$.
\ere

But \dd\dpdb generic pairs are not necessarily
generic in the sense of the true forcing product
$\dP\ti\dP$.
Indeed, if say $\dP=\text{Sacks}$ (all perfect trees)
then any \dd\dpd generic pair $\ang{x,y}$ has the property
that $x,y$ belong to same \dd\Eo invariant Borel sets
coded in the ground universe, while for any uncountable
and co-uncountable Borel set $U$ coded in the ground
universe there is a \dd{\dP\ti\dP}generic pair
$\ang{x,y}$ with $x\in U$ and $y\nin U$.

\ble
\lam{stf2+}
Assume that\/ $\dP$ is a\/ \sif, $n\ge1$,
$P\in\sct\dP n$, and a set\/ $D\sq\dPd$ is open dense
in\/ $\dPd$.
Then there is a tree\/ $Q\in\sct\dP n$
such that\/ $Q\nq n P$ and\/ $\ang{\raw{Q}{s},\raw Qt}\in D$
\iman
whenever\/ $s,t\in 2^n$ and\/ $s(n-1)\ne t(n-1)$.
\ele
\bpf
[compare to Lemma~\ref{stf}\ref{stf2c}]
Let $s,t\in 2^n$ be any pair with $s(n-1)\ne t(n-1)$.
By the density there is a condition $\ang{S,T}\in D$
such that $S\sq \raw{P}{s}$ and $T\sq\raw Pt$.
Note that $T=\sg\app S$ for some $s\in\bse$ since
$\ang{S,T}\in\dpd$.
Applying Lemma~\ref{stf}\ref{stf3}
($n+1$ there corresponds to $n$ here)
we obtain a tree $P'\in\sct\dP n$
such that $P'\nq n P$ and
$\raw{P'}{s}\sq S$, $\raw{P'}t\sq T$.
Then $\ang{\raw{P'}{s},\raw{P'}t}\in D$, as $D$ is open.
Consider all pairs $s,t\in 2^n$ with $s(n-1)\ne t(n-1)$
one by one.
\epf

\ble
\lam{dr}
Assume that\/ $\dP$ is a\/ \sif, $\ang{T,T'}\in\dpd$, $n<\om$,
$s,t\in2^n$.
Then\/ $\ang{\raw Ts,\raw {T'}t}\in\dpd$.
\ele
\bpf
Let $\sg\in\bse$ satisfy $\sg\app T=T'$.
Note that $\oi T=\oi {T'}$, hence we define
$h=\oin Tn=\oin{T'}n$.
By Lemma~\ref{sadd}\ref{sadd3}, there are strings
$u\in 2^h\cap T$ and $v\in 2^h\cap T'$ such that
$\raw Ts=\req Tu$ and $\raw{T'}t=\req{T'}v$.
Then obviously $\sg\app \req Tu=\req {T'}{v'}$,
where $v'=\sg\app u$.
On the other hand
$\req{T'}v=(v\app v')\app \req{T'}{v'}$ by
Lemma~\ref{trans}.
It follows that $\req{T'}v=(v\app v'\app\sg)\app \req Tu$, as
required.
\epf

\bcor
\lam{lnr}
Assume that\/ $\dP$ is a\/ \sif.
Then\/ $\dpd$ forces\/ $\rpi\lef \not\Eo \rpi\rig$, where
\iman
$\ang{\rpi\lef,\rpi\rig}$ is a name of the\/ \dd\dpdb generic
pair.
\ecor
\bpf
Otherwise a condition $\ang{T,T'}\in\dpd$ forces
\iman
$\rpi\rig=\sg\app\rpi\lef$, where $\sg\in\bse.$
Find $n$ and $s,t\in2^n$ such that
$\raw {T'}t\cap (\sg\app \raw Ts)=\pu$
and apply the lemma.
\epf

\section{\Mut s}
\las{mut}

Let a \rit{\mut} be any sequence
\index{multitree}%
$\vpi=\sis{\ang{\mtt\vpi k,\mtp\vpi k}}{k<\om}$ such that\vom
\index{zztauphik@$\mtt\vpi k$}%
\index{zzpphik@$\mtp\vpi k$}%
\ben
\snenu
\itla{mut1}
if $k<\om$ then $\mtp\vpi k\in\om\cup\ans{-1}$, and
\imar{mut1}
the set $\abs\vpi=\ens{k}{\mtp\vpi k\ne-1}$
\index{zz:phi:@$\abs\vpi$}
(the \rit{support} of $\vpi$) is finite;\vom
\index{support!@$\abs\vpi$}%


\itla{mut2}
if $k\in\abs\vpi$ then
\imar{mut2}
$\mtt\vpi k=
\ang{\ntt\vpi k0,\ntt\vpi k1,\dots,\ntt\vpi k{\mtp\vpi k}}$,
\index{zzTphikn@$\ntt\vpi kn$}
where each $\ntt\vpi kn$ is a tree in $\pes$ and
$\ntt\vpi k{n}\nq n \ntt\vpi k{n-1}$
whenever $1\le n\le \mtp\vpi k$,
while if $k\nin\abs\vpi$ then simply $\mtt\vpi k=\La$
(the empty sequence).\vom
\een
In this context, if $n \le \mtp\vpi k$ and $s\in2^n$ then
let $\ntt\vpi ks=\raw{\ntt\vpi kn}{s}$.
\index{zzTphiks@$\ntt\vpi ks$}

Let $\vpi,\psi$ be \mut s.
Say that
$\vpi$ \rit{extends} $\psi$, symbolically $\psi\cle\vpi$, if
$\abs \psi\sq\abs\vpi$, and, for every $k\in\abs\psi$, we have
\index{multitree!extends}%
\index{multitree!$\psi\cle\vpi$}%
\index{zzpsi<vpi@$\psi\cle\vpi$}%
$\mtp\vpi k\ge \mtp\psi k$ and
$\mtt\vpi k$ extends $\mtt\psi k$, so that
$\ntt\vpi k{n}=\ntt\psi k{n}$ for all $n\le \mtp\psi k$;

If $\dP$ is a \sif\ then let $\mt\dP$ (\rit{\mut s over $\dP$})
\index{zzMTP@$\mt\dP$}%
be the set of all \mut s
$\vpi$ such that $\ntt\vpi kn\in\sct\dP n$
whenever $k\in\abs\vpi$ and $n\le \mtp\vpi k$.

\section{Jensen's extension of a large-tree forcing notion}
\las{jex}

Let $\zfc'$ be the subtheory of $\zfc$ including all
\index{zzzfc@$\zfc'$}%
axioms except for the
power set axiom, plus the axiom saying that $\pws\om$ exists.
(Then $\omi$, $\dn$, and sets like $\pet$ exist as well.)

\bdf
\lam{dPhi}
Let $\cM$ be a countable transitive model of $\ZFC'$.
\index{zzMgot@$\cM$}%
\index{model!Mgot@$\cM$}%
Suppose that $\dP\in\cM\yd \dP\sq\pes$ is a \sif.
Then $\mt\dP\in\cM$.
A set $D\sq\mt\dP$ is \rit{dense in $\mt\dP$} iff
\index{set!dense}%
for any $\psi\in\mt\dP$ there is
a \mut\ $\vpi\in D$ such that $\psi\cle\vpi$.

Consider any \dd\cle increasing sequence
$\dphi=\sis{\vpi(j)}{j<\om}$ of \mut s
\index{zzPhit@$\dphi$}%
$$
\vpi(j)=
\sis{\ang{\mtt{\vpi(j)} k,\mtp{\vpi(j)} k}}{k<\om}\in\mt\dP\,,
$$
\rit{generic over\/ $\cM$} in the sense that it intersects
every set $D$, $D\sq\mt\dP$, dense in $\mt\dP$,
which belongs to $\cM$.
Then in particular $\dphi$ intersects every set
$$
D_{kp}=\ens{\vpi\in\mt\dP}
{k\in\abs{\vpi}\land\mtp\vpi k\ge p}\,,\quad k,p<\om\,.
$$
Therefore if $k<\om$ then by definition
there is an infinite sequence
$$
\dots \nq 5 \tx k{4}\nq 4 \tx k{3}
\nq 3 \tx k{2}\nq 2 \tx k{1}\nq 1 \tx k{0}
$$
of trees $\tx k{n} \in\sct\dP n$, such that, for any $j$, if
\index{zzTPhikn@$\tx k{n}$}%
$k \in\abs{\vpi(j)}$ and $n\le\mtp{\vpi(j)}k$ then
$\ntt{\vpi(j)} kn = \tx k{n}$.
If $n<\om$ and $s\in2^n$ then we let
$\tx ks=\raw{\tx kn}s$;
\index{zzTPhiks@$\tx k{s}$}%
then $\tx ks\in \dP$ since $\tx kn\in\sct\dP n$.
Then it follows from Lemma~\ref{fus} that
$$
\TS
\ufi k=\bigcap_{n}\tx kn=\bigcap_{n}\bigcup_{s\in2^n}\tx ks
\index{zzUPhik@$\ufi k$}%
\eqno(1)
$$
is a tree in $\pes$ (not necessarily in $\dP$), as well as
the trees $\uf ks$,
and still by Lemma~\ref{fus},
\index{zzUPhiks@$\uf k{s}$}%
$$
\TS
\uf ks=\ufi k\cap \tx ks
=\bigcap_{n\ge \lh s}\raw{\tx kn}s
=\bigcap_{n\ge \lh s}\bigcup_{t\in2^n,\:s\sq t}\tx kt \,,
\eqno(2)
$$
and obviously $\ufi k=\uf k\La$.

Define a set of trees
$\dU=\ens{\sg\app \uf ks}{k<\om\land s\in\bse\land
\sg\in\bse}\sq\pes$.%
\index{zzU@$\dU$}%
\edf

The next few simple lemmas show useful effects of the
genericity of $\dphi$; their common motto is that the
extension from $\dP$ to $\dP\cup\dU$ is rather innocuous.

\ble
\lam{uu1}
Both\/ $\dU$ and the union\/
$\dP\cup\dU$ are\/ \sif s$;$ $\dP\cap\dU=\pu$.
\ele
\bpf
To prove the last claim, let $T\in\dP$ and $U=\uf ks\in \dU$.
(If $U=\sg\app\uf ks$, $\sg\in\bse,$ then replace
$T$ by $\sg\app T$.)
The set $D(T,k)$ of all \mut s $\vpi\in \mt\dP$,
such that $k\in\abs\vpi$ and $T\bez\raw{\ntt\vpi kn}s\ne\pu$, where
$n=\mtp\vpi k$, belongs to $\cM$
and obviously is dense in $\mt\dP$.
Now any \mut\ $\vpi(j)\in D(T,k)$ witnesses that
$T\bez\uf ks\ne\pu$.
\epf

\ble
\lam{uu2}
The set\/ $\dU$ is dense in\/ $\dU\cup\dP$.
The set\/ $\dud$ is dense in\/\break $\dpud$.
\ele
\bpf
Suppose that $T\in\dP$.
The set $D(T)$ of all \mut s $\vpi\in \mt\dP$,
such that $\ntt\vpi k0=T$ for some $k$,
belongs to $\cM$
and obviously is dense in $\mt\dP$.
It follows that $\vpi(j)\in D(T)$ for some $j$,
by the choice of $\dphi$.
Then $\tx k\La=T$ for some $k$.
However by construction $\uf k\La=\ufi k\sq \tx k\La$.

Now suppose that $\ang{T,T'}\in\dPd$, so that $T'=\sg\app T$,
$\sg\in\bse$.
By Lemma~\ref{uu1} ($\dP\cap\dU=\pu$) it is impossible that
one of the trees $T,T'$ belongs to $\dP$ and the other one
to $\dU$.
Therefore we can assume that $T,T'\in\dP$.
By the first claim of the lemma, there is a tree $U\in\dU$,
$U\sq T$.
Then $U'=\sg\app U\in\dU$ and still $U'=\sg\app U$, hence
$\ang{U,U'}\in\dUd$, and it extends $\ang{T, T'}$.
\epf


\ble
\lam{uu22}
If\/ $k,l<\om$, $k\ne l$, and\/ $\sg\in\bse$ then\/
$\ufi k\cap (\sg\app\ufi l)=\pu$.
\ele
\bpf
The set $D'(k,l)$ of all \mut s $\vpi\in \mt\dP$,
such that $k,l\in\abs\vpi$ and
$\ntt\vpi kn\cap (\sg\app\ntt\vpi lm)=\pu$
for some $n\le \mtp\vpi k$, $m\le \mtp\vpi l$,
belongs to $\cM$ and is dense in $\mt\dP$.
So $\vpi(j)\in D'(k,l)$ for some $j<\omega.$
But then for some $n,m$ we have
$\ufi k\cap (\sg\app\ufi l)
\subseteq \ntt{\vpi(j)} kn\cap (\sg\app\ntt{\vpi(j)} lm)=\pu$.
\epf

\bcor
\lam{uu2c}
\imac
If\/ $\ang{U,U'}\in\dud$ then there exist$:$ $k<\om$,
strings\/ $s,s'\in\bse$ with\/ $\lh s=\lh{s'}$,
and strings\/ $\sg,\sg'\in\bse,$ such that\/ $U=\sg\app\uf ks$
and\/ $U'=\sg'\app\uf k{s'}$.
\ecor
\bpf
By definition, we have $U=\sg\app\uf ks$
and\/ $U'=\sg'\app\uf {k'}{s'}$, for suitable $k,k'<\om$
and $s,s',\sg,\sg'\in\bse.$
As $\ang{U,U'}\in\dud$, it follows from Lemma~\ref{uu22}
that $k'=k$, hence $U'=\sg\app\uf {k}{s'}$.
Therefore $\sg\app\uf ks=\tau\app\sg'\app\uf {k}{s'}$
for some $\tau\in\bse.$
In other words, $\uf ks=\tau'\app\uf {k}{s'}$, where
$\tau'=\sg\app\sg'\app\tau\in\bse.$
It easily follows that $\lh s=\lh{s'}$.
\epf

The two following lemmas  show that, due to
the generic character of extension, those pre-dense sets
which belong to $\cM$, remain pre-dense in the extended
forcing.

Let $X\sqf\bigcup D$ mean that there is a finite set
$D'\sq D$ with $X\sq\bigcup D'$.

\ble
\lam{uu3}
If a set\/ $D\in\cM$, $D\sq\dP$ is pre-dense
in\/ $\dP$, and\/ $U\in\dU$, then\/ $U\sqf\bigcup D$.
\imar{\color{red}split\\ into 2 lemmas}%
Moreover\/ $D$ is pre-dense in\/ $\dU\cup\dP$.
\ele
\bpf
We can assume that $D$ is in fact open dense in $\dP$.
(Otherwise replace it with the set
$D'=\ens{T\in\dP}{\sus S\in D\,(T\sq S)}$ which also
belongs to $\cM$.)

We can also assume that $U=\uf ks\in\dU$,
where $k<\om$ and $s\in\bse.$
(The general case,
when $U=\sg\app\uf ks$ for some $\sg\in\bse,$
is reducible to the case $U=\uf ks$ by
substituting the set $\sg\app D$ for $D$.)
\vyk{
To prove the first claim,
suppose that $U=\ufi k\in\dU$, $k<\om$.
(The general case,
when $U=\sg\app\uf ks$ for some $s,\sg\in\bse,$
is reducible to the intermediate case $U=\uf ks$ by
substituting the set $\sg\app D$ for $D$;
then we simply drop $s$
since $\uf ks\sq\uf k\La=\ufi k$.)
}%

The set $\Da\in\cM$ of all \mut s
$\vpi\in \mt\dP$ such that
$k\in\abs\vpi$,
$\lh s<h=\mtp{\vpi}k$, and
$\raw{\ntt\vpi kh}t\in D$ for all $t\in 2^{h}$,
is dense in $\mt\dP$ by Lemma~\ref{stf}\ref{stf2c} and
the open density of $D$.
Therefore there is an index $j$ such that
$\vpi(j)\in\Da$.
Let $h(j)=\mtp{\vpi(j)}k$.
Then the tree
$S_t=\raw{\ntt{\vpi(j)} k{h(j)}}t=\raw{\tx k{h(j)}}t=\tx kt$
belongs to $D$
for all $t\in 2^{h(j)}$.
We conclude that
$$
\textstyle
U=\uf ks
\sq\ufi k
\sq\bigcup_{t\in 2^{h(j)}}\tx kt
\sq\bigcup_{t\in 2^{h(j)}}S_t=\bigcup D'\,,
$$
where $D'=\ens{S_t}{t\in 2^{h(j)}}\sq D$
is finite.

To prove the pre-density claim,
\vyk{
consider any tree $U=\uf ks\in\dU$, where $k<\om$, $s\in\bse.$
(The general case,
when $U=\sg\app\uf ks$ for some $\sg\in\bse,$
is reducible to the case $U=\uf ks$ by
substituting $\sg\app D$ for $D$.)
Let $j$ and $h(j)$ be defined as above.
P
}%
pick a string $t\in 2^{h(j)}$ with $s\su t$.
Then $V=\uf kt\in\dU$ and $V\sq U$.
However $V\sq \tx kt=S_t\in D$.
Thus  $V$ witnesses that $U$ is compatible with
$S_t\in D$ in $\dU\cup\dP$, as required.
\epf

\ble
\lam{uu3d}
If a set\/ $D\in\cM$, $D\sq\dPd$ is pre-dense
in\/ $\dPd$ then\/ $D$ is pre-dense in\/ $\dpud$.
\iman
\ele
\bpf
Let $\ang{U,U'}\in\dud$; the goal is to prove that
$\ang{U,U'}$ is compatible in $\dpud$ with a condition
$\ang{T,T'}\in D$.
By Corollary~\ref{uu2c}, there exist: $k<\om$ and strings
$s,s',\sg,\sg'\in\bse$ such that $\lh s=\lh{s'}$ and
$U=\sg\app\uf ks$, $U'=\sg'\app\uf k{s'}$.
As in the proof of the previous lemma, we can assume
that $\sg=\sg'=\La$, so that $U=\uf ks$, $U'=\uf k{s'}$.
(The general case is reducible to this case by substituting
the set
$\ens{\ang{\sg\app T,\sg'\app T'}}{\ang{T,T'}\in D}$
for $D$.)

Assume that $D$ is in fact open dense.

Consider the set $\Da\in\cM$ of all \mut s
$\vpi\in \mt\dP$ such that
$k\in\abs\vpi$, $\lh s=\lh{s'}=n< h=\mtp{\vpi}k$,
and $\ang{\raw{\ntt\vpi kh}u,\raw{\ntt\vpi kh}{u'}} \in D$
whenever $u,u'\in 2^{h}$ and $u(h-1)\ne u'(h-1)$.
The set $\Da$
is dense in $\mt\dP$ by Lemma~\ref{stf2+}.
Therefore $\vpi(j)\in\Da$ for some $j$, so that
if $u,u'\in 2^{h(j)}$, where $h(j)=\mtp{\vpi(j)}k> n$,
and $u({h(j)}-1)\ne u'({h(j)}-1)$, then
$$
\ang{\raw{\ntt{\vpi(j)}k{h(j)}}u,\raw{\ntt{\vpi(j)}k{h(j)}}{u'}}
=\ang{\tx ku,\tx k{u'}}\in D\,.
$$
Now, as ${h(j)}>n$, let us pick $u,u'\in 2^{h(j)}$ such that
$u({h(j)}-1)\ne u'({h(j)}-1)$ and $s\su u$, $s'\su u'$.
Then $\ang{\tx ku,\tx k{u'}}\in D$.
On the other hand, the pair
$\ang{\uf ku,\uf k{u'}}$ belongs to $\dud$ by
Lemma~\ref{dr},
$$
\ang{\uf ku,\uf k{u'}}\le \ang{\uf ks,\uf k{s'}}\,,
$$
and
finally we have
$\ang{\uf ku,\uf k{u'}}\le \ang{\tx ku,\tx k{u'}}$.
We conclude that the given condition
$\ang{\uf ks,\uf k{s'}}$
is compatible with the condition
$\ang{\tx ku,\tx k{u'}}\in D$, as required.
\epf

\section{Real names}
\las{rena}

In this Section, we
assume that $\dP$ is a \sif\ and $\bse\in\dP$.
It follows by \ref{sitf2}
that all trees $T[s]=\raw{(\bse)}s$ (see Example~\ref{clop})
also belong to $\dP$.

Recall that $\dPd$ adds a pair of reals
$\ang{x\lef,x\rig}\in\dn\ti\dn$.

Arguing in the conditions of Definition~\ref{dPhi},
the goal of the following Theorem~\ref{K} will be to prove that,
for any \dd\dpdb name $c$ of a real in $\dn,$
it is forced by the extended
forcing $\dpud$ that $c$ does not belong to sets
of the form $[U]$, where
$U$ is a tree in $\dU$, {\ubf unless} $c$ is a name of one of
reals in the \dd\Eo class of one of the generic reals
$x\lef,x\rig$ themselves.

We begin with a suitable notation.

\bdf
\lam{rk}
A \rit{\dd\dpdb real name} is a system
\index{real name}%
\index{PxPreal name@\dd\dpdb real name}%
$\rc=\sis{\kc ni}{n<\om,\, i<2}$ of sets $\kc ni\sq\dpd$
such that each set $C_n=\kc n0\cup \kc n1$ is
pre-dense in $\dPd$
and any conditions $\ang{S,S'}\in \kc n0$ and
$\ang{T,T'}\in \kc n1$ are
incompatible in $\dPd$.

If a set $G\sq\dPd$ is \dd\dpdb generic at least over
the collection of all  sets $C_n$ then we define
$\rc[G]\in\dn$ so that $\rc[G](n)=i$ iff $G\cap \kc ni\ne\pu$.
\index{zzc:G:@$\rc[G]$}%
\edf

Any \dd\dpdb real name $\rc=\sis{\kc ni}{}$
induces (can be understood as) a \dd\dpdb name
(in the ordinary forcing notation)
for a real in $\dn.$

\bdf
[actions]
\lam{rka}
Strings in $\bse$ can act on names
$\rc=\sis{\kc ni}{n<\om\yi i<2}$
in two ways,
related either to conditions or to the output.

If $\sg,\sg'\in\bse$
then define a \dd\dpdb real name
$\ang{\sg,\sg'}\apc\rc=\sis{\ang{\sg,\sg'}\app\kc ni}{}$,
\index{zzsgsg.c@$\ang{\sg,\sg'}\apc\rc$}%
\index{action!sgsg.c@$\ang{\sg,\sg'}\apc\rc$}%
where
$\ang{\sg,\sg'}\app\kc ni=
\ens{\ang{\sg\app T,\sg'\app T'}}{\ang{T,T'}\in\kc ni}$
for all $n,i$.

If $\rho\in\bse$ then define a \dd\dpdb real name
$\rho\app\rc= \sis{\rkc ni\rho}{}$,
\index{zzrho.c@$\rho\app\rc$}%
\index{action!rho.c@$\rho\app\rc$}%
where
$\rkc ni\rho=\kc n{1-i}$ whenever $n<\lh\rho$ and $\rho(n)=1$,
but $\rkc ni\rho=\kc n{i}$ otherwise.
\edf

Both actions are idempotent.
The difference between them is as follows.
If $G\sq\dPd$ is a \dd\dpdb generic set then
$(\ang{\sg,\sg'}\apc\rc)[G]=\rc[\ang{\sg,\sg'}\apc G]$,
where
$\ang{\sg,\sg'}\apc G=
\ens{\ang{\sg\app T,\sg'\app T'}}{\ang{T,T'}\in G}$,
while $(\rho\app\rc)[G]=\rho\app (\rc[G])$.

\bpri
\lam{proj}
Define a \dd\dpdb real name $\rpi\lef=\sis{\kc ni}{n<\om\yi i<2}$
\index{zzxleft@$\rpi\lef$}%
\index{zzxright@$\rpi\rig$}%
such that each set $\kc ni\sq\dpd$ contains all pairs of the
form $\ang{T[s],T[t]}$, where $s,t\in2^{n+1}$ and $s(n)=i$,
and a \dd\dpdb real name
$\rpi\rig=\sis{\kc ni}{n<\om\yi i<2}$
such that accordingly each
set $\kc ni\sq\dpd$ contains all pairs
$\ang{T[s],T[t]}$, where $s,t\in2^{n+1}$ and now $t(n)=i$.
%
\epri

Then $\rpi\lef$, $\rpi\rig$ are names of the
\dd\dP generic reals $x\lef$, resp., $x\rig$,
and each name $\sg\app\rpi\lef$
($\sg\in\bse$) induces a \dd\dpdb name of the real
$\sg\app (x\lef[G])$;
the same for ${}\rig$.

\section{Direct forcing a real to avoid a tree}
\las{saway2}

Let $\rc=\sis{\kc ni}{}$, $\rd=\sis{\kd ni}{}$
be \dd\dpdb real names.
Say that a condition $\ang{T,T'}\in \pes\ti_{\Eo}\pes$:\vom
\bit
\item
\rit{directly forces\/ $\rc(n)=i$},
where $n<\om$, $i=0,1$, if $\ang{T,T'}\le\ang{S,S'}$
for some $\ang{S, S'}\in\kc ni$;\vom

\item
\rit{directly forces\/ $s\su\rc$},
where $s\in\bse,$ iff for all $n<\lh s$, $\ang{T,T'}$
directly forces $\rc(n)=i$, where $i=s(n)$;\vom

\item
\rit{directly forces\/ $\rd\ne\rc$}, iff there are strings
$s,t\in\bse,$ incomparable in $\bse$ and such that
$\ang{T,T'}$ directly forces $s\su\rc$ and $t\su\rd$;\vom

\item
\rit{directly forces\/ $\rc\nin[U]$},
where $U\in\pet$, iff there is a string $s\in\bse\bez U$
such that $\ang{T,T'}$ directly forces $s\su \rc$.
\eit


\ble
\lam{K1}
If\/ $S\in\dP$, $\ang{R,R'}\in\dpd$, and\/
$\rc$ is a\/ \dd\dpdb real name, then there exists a tree\/
$S'\in\dP$ and a condition\/
$\ang{T,T'}\in\dpd$, $\ang{T,T'}\leq\ang{R,R'}$,
such that\/ $S'\sq S$ and\/
$\ang{T,T'}$ directly forces\/ $\rc\nin[S']$.
\ele
\bpf
Clearly there is a condition\/
$\ang{T,T'}\in\dpd$, $\ang{T,T'}\leq\ang{R,R'}$, which
directly forces\/ $u\su\rc$ for some $u\in\bse$
satisfying $\lh u>\lh{(\roo S)}$.
There is a string  $v\in S$, $\lh v=\lh u$,
incomparable with $u$.
The tree $S'=\req S{v}$ belongs to $\dP$, $S'\sq S$
by construction, and obviously
$\ang{T,T'}$ directly forces $\rc\nin [S']$.
\epf

\ble
\lam{K2}
If\/ 
$\rc$ is a\/ \dd\dpdb real name, $\sg\in\bse,$
and a condition\/ $\ang{R,R'}\in\dpd$ directly forces\/
$\sg\app\rc\ne\rpl$, resp., $\sg\app\rc\ne\rpr$,
then there is a stronger condition\/
$\ang{T,T'}\in\dpd$, $\ang{T,T'}\leq\ang{R,R'}$,
which directly forces resp.\ $\rc\nin[\sg\app T],$
$\rc\nin[\sg\app T'].$
\ele
\bpf
We just prove the ``left'' version,
as the ``right'' version can be proved similarly.
So let's assume  that
$\ang{R,R'}$ directly forces $\rc\ne\rpl$.
There are incomparable strings $u,v\in\bse$
such that $\ang{R,R'}$ directly forces
$u\su\sg\app\rc$, hence, $\sg\app u\su\rc$ as well, and
also directly forces $v\su\rpl$.
Then by necessity $v\in R$, hence $T=R\ret v\in \dP$,
but $u\nin T$.
Let $T'=\rho\app T$, where $\rho\in\bse$ satisfies
$R'=\rho\app R$.
By definition, the condition
$\ang{T,T'}\in\dpd$ directly forces\/ $\rc\nin[\sg\app T]$
(witnessed by $s=\sg\app u$),
as required.
\epf

\bte
\lam{K}
With the assumptions of Definition~\ref{dPhi}, suppose that\/
$\rc=\break
\sis{\kc mi}{m<\om\yi i<2}\in\cM$ is a\/ \dd\dpdb real name,
and for every\/ $\sg\in\bse$ the set
$$
D_\sg=
\ens{\ang{T,T'}\in\dPd}
{\ang{T,T'}\text{ directly forces\/ }\rc\ne{\sg\app\rpl}
\text{ and\/ }\rc\ne{\sg\app\rpr}}
$$
is dense in\/ $\dPd$.
Let\/ $\ang{W,W'}\in\dpud$ and\/
$U\in \dU$.

Then there is a stronger condition\/
$\ang{V,V'}\in\dUd\yd \ang{V,V'}\le \ang{W,W'}$, which
directly forces\/ $\rc\nin[U]$.
\ete
\bpf
By construction, $U=\rho\app\uf K{s_0}$,
where $K<\om$ and $\rho,s_0\in\bse;$
we can assume that simply $s_0=\La$, so that $U=\rho\app\ufi K$.
Moreover
we can assume that $\rho=\La$ as well, so that $U=\ufi K$
(for if not then replace $\rc$ with $\rho\app\rc$).

Further, by Corollary~\ref{uu2c}, we can assume that
$W=\sg\app\uf L{t_0}\in\dU$ and
$W'=\sg'\app\uf L{t_0'}\in\dU$, where $L<\om$,
$t_0,t_0'\in\bse,$ $\lh {t_0}=\lh {t_0'}$,
and $\sg,\sg'\in\bse.$
And moreover we can assume that $\sg=\sg'=\La$, so that
$W=\uf L{t_0}$ and $W'=\uf L{t_0'}$
(for if not then replace $\rc$ with $\ang{\sg,\sg'}\apc\rc$).

The indices $K,L$ involved can be either equal or different.

There is an index $J$ such that the \mut\
$\vpi(J)$ satisfies $K,L\in\abs{\vpi(J)}$ and
$\mtp{\vpi(J)}L\ge{h_0}=\lh{t_0}=\lh{t'_0}$,
so that the trees $S_0=\ntt{\vpi(J)}K0=\tx K{0}$,
$$
T_0=\raw{\ntt{\vpi(J)}L{h_0}} {t_0}=\tx L{t_0}\,,\quad
T'_0=\raw{\ntt{\vpi(J)}L{h_0}} {t_0'}=\tx L{t_0'}
$$
in $\dP$ are defined.
Note that $U\sq S_0$ and $W\sq T_0$, $W'\sq T_0'$
under the above assumptions.

Let $\cD$ be the set of all \mut s
$\vpi\in \mt\dP$ such that
$\vpi(J)\cle\vpi$ and for every pair $t,t'\in 2^{n}$,
where $n=\mtp{\vpi}{L}$, such that $t(n-1)\ne t'(n-1)$,
the condition $\ang{\ntt\vpi L{t},\ntt\vpi L{t'}}$
directly forces $\rc\nin[\ntt{\vpi}Km]$, where
$m=\mtp{\vpi}{K}$.
\iman

\bkl
\lam{Ka}
$\cD$ is dense in\/ $\mt\dP$ above\/ $\vpi(J)$.
\ekl
\bpf
Let a \mut\ $\psi\in\mt\dP$ satisfy
\iman
$\vpi(J)\cle\psi$; the goal is to define a \mut\
$\vpi\in\cD$, $\psi\cle\vpi$.
Let $m=\mtp{\psi}{K}$, $n=\mtp{\psi}{L}$,
$Q=\ntt{\psi}K{m}$, $P=\ntt{\psi}L{n}$.
\vyk{
Pick strings $r,r'\in 2^{n}$ with $t_0\su r$, $t'_0\su r'$,
and let
$$
R=\ntt{\psi}L{r\we0}=\raw P{r\we0}\,,\quad
R'=\ntt{\psi}L{r'\we1}=\raw P{r'\we1}\,.
$$
Then
$R,R'\sq T_0$ are trees in $\dP$, $\ang{R,R'}\in\dpd$
(Lemma~\ref{dr}),
$\ang{R,R'}\leq\ang{T_0,T_0'}$.
}%
\vom

{\it Case 1\/}: $K\ne L$.
Consider any $s\in 2^{m+1}$ and $t,t'\in 2^{n+1}$
with $t(n)\ne t'(n)$.
By Lemma~\ref{K1}, there is a tree $S\in \dP$
and a condition $\ang{R,R'}\in\dpd$ such that
$S\sq \raw Qs$, $\ang{R,R'}\le\ang{\raw Pt,\raw P{t'}}$,
and $\ang{R,R'}$ directly forces\/ $\rc\nin[S]$.
By Lemma~\ref{stf}\ref{stf2},\ref{stf3} there are trees
$Q_1\in\sct \dP{m+1}$ and $P_1\in\sct \dP{n+1}$
such that $Q_1\nq{m+1}Q$, $P_1\nq{n+1}P$,
$\raw{Q_1}s=S$ and
$\ang{\raw{P_1}t,\raw{P_1}{t'}}\le \ang{R,R'}$.

Repeat this procedure
so that all strings
$s\in 2^{m+1}$ and all pairs of strings
$t,t'\in 2^{n+1}$ with $t(n)\ne t'(n)$ are considered.
We obtain trees $Q'\in\sct \dP{m+1}$ and $P'\in\sct \dP{n+1}$
such that $Q'\nq{m+1}Q$, $P'\nq{n+1}P$, and if
$s\in 2^{m+1}$ and $t,t'\in 2^{n+1}$, $t(n)\ne t'(n)$, the
condition $\ang{\raw{P'}t,\raw{P'}{t'}}$
directly forces $\rc\nin[\raw{Q'}{s}]$ --- hence
directly forces $\rc\nin[Q']$.

Now define a \mut\ $\vpi\in\mt\dP$ so that
$\abs{\vpi}=\abs\psi$,
$\mtp{\vpi}k=\mtp{\psi}k$ and $\mtt{\vpi}k=\mtt{\psi}k$
for all $k\nin\ans{K,L}$,
$\mtp{\vpi}{K}=m+1$,
$\mtp{\vpi}{L}=n+1$,
and
$\ntt\vpi K{m+1}=P'$, $\ntt\vpi L{n+1}=Q'$ as the new elements
of the $K$th and $L$th components.
We have $\vpi\in \cD$ and $\psi\cle\vpi$ by construction.
(Use the fact that $P'\nq{n+1}P$ and $Q'\nq{m+1}Q$.)
\vom

{\it Case 2\/}: $L=K$, and hence $m=n$ and $P=Q$.
Let $h=\oin Pn$.
Consider any pair $t,t'\in 2^{n+1}$ with $t(n)\ne t'(n)$.
In our assumptions there is a condition
$\ang{U,U'}\in\dpd$, $\ang{U,U'}\leq\ang{\raw{T}t,\raw{T}{t'}}$,
which directly forces
both $\rc\ne{\sg\app\rpl}$ and $\rc\ne{\sg\app\rpr}$
for any $\sg\in 2^h$.
By Lemma~\ref{K2}, there is a stronger condition
$\ang{T,T'}\in\dpd$, $\ang{T,T'}\leq\ang{U,U'}$, which
directly forces both $\rc\nin[\sg\app T]$ and
$\rc\nin[\sg\app T']$ still for all $\sg\in 2^h$.
Then as in Case 1, there is a tree
$P_1\in\sct\dP{n+1}$, $P_1\nq{n+1}P$, such that
$\raw{P_1}{t}\sq T$, $\raw{P_1}{t'}\sq T'$.

We claim that $\ang{T,T'}$ directly forces
$\rc\nin[P_1]$, or equivalently, directly forces
$\rc\nin[\raw{P_1}{s\we i}]$ for any $s\we i\in 2^{n+1}$
(then $s\in 2^n$).
Indeed if $s\we i\in 2^{n+1}$ then
$\raw{P_1}{s\we i} =\sg\app \raw{P_1}{t}$ or
$=\sg\app \raw{P_1}{t'}$ for
some $\sg\in2^h$ by the choice of $h$.
Therefore $\raw{P_1}{s\we i}$ is a subtree of one of
the two trees $\sg\app T$ and $\sg\app T'$.
The claim now follows from the choice of $\ang{T,T'}$.
We conclude that the stronger condition
$\ang{\raw{P_1}{t},\raw{P_1}{t'}}\le \ang{T,T'}$ also
directly forces $\rc\nin[P_1]$.

Repeat this procedure so that all pairs of strings
$t,t'\in 2^{n+1}$ with $t(n)\ne t'(n)$ are considered.
We obtain a tree $P'\in\sct \dP{n+1}$
such that $P'\nq{n+1}P$, and if
$t,t'\in 2^{n+1}$, $t(n)\ne t'(n)$, then
$\ang{\raw{P'}t,\raw{P'}{t'}}$
directly forces $\rc\nin[P']$.

Similar to Case 1, define a \mut\ $\vpi\in\mt\dP$ so that
$\abs{\vpi}=\abs\psi$,
$\mtp{\vpi}k=\mtp{\psi}k$ and $\mtt{\vpi}k=\mtt{\psi}k$
for all $k\ne K$,
$\mtp{\vpi}{K}=n+1$, and $\ntt\vpi K{n+1}=P'$ as the new
element of the $(K=L)$th component.
Then $\vpi\in \cD$, $\psi\cle\vpi$.
\epF{Claim}

We come back to the proof of Theorem~\ref{K}.
The lemma implies that there is an index $j\ge J$
such that the \mut\ $\vpi(j)$ belongs to $\cD$.
Let $n=\mtp{\vpi(j)}{L}$, $m=\mtp{\vpi(j)}{K}$.
Pick strings $t,t'\in 2^{n}$ such that $t_0\su t$,
$t_0'\su t'$, $t(n)\ne t'(n)$.
Let
\bce
$T=\ntt{\vpi(j)}L{t}=\tx Lt$,
$T'=\ntt{\vpi(j)}L{t'}=\tx L{t'}$,
$S=\ntt{\vpi(j)}K{m}=\tx Km$.
\ece
Then $\ang{T,T'}\in\dPd\yt \ang{T,T'}\le \ang{T_0,T'_0}$,
and
$\ang{T,T'}$ directly forces $\rc\nin[S]$.

Consider the condition $\ang{V,V'}\in\dud$, where
$V=\uf L{t}$ and $V'=\uf L{t'}$ belong to $\dU$.
(Recall that $V=\raw{\ufi L}t$ and $V'=\raw{\ufi L}{t'}$,
and hence $V'=\sg\app V$ for a suitable $\sg\in\bse$.)
By construction we have both $\ang{V,V'}\leq \ang{W,W'}$
(as $t_0\sq t,t'$)
and $\ang{V,V'}\leq \ang{T,T'}\leq\ang{T_0,T'_0}$.
Therefore $\ang{V,V'}$ directly forces $\rc\nin[S]$.
And finally, we have $U\sq \ntt{\vpi(j)}Km=S$, so that
$\ang{V,V'}$ directly forces $\rc\nin[U]$, as required.
\epF{Theorem~\ref{K}}

\section{Jensen's forcing}
\las{jfor}

In this section,
{\ubf we argue in $\rL$, the constructible universe.}
Let $\lel$ be the canonical wellordering of $\rL$.

\bdf
[in $\rL$]
\lam{uxi}
Following the construction in
\cite[Section 3]{jenmin} \rit{mutatis mutandis},
define, by induction on $\xi<\omi$, a countable \sif\
$\dU_\xi\sq\pes$
as follows.

Let $\dU_0$ consist of all trees of the form $T[s]$,
see Example~\ref{clop}.

Suppose that $0<\la<\omi$, and countable \sif s
$\dU_\xi\sq\pes$ are defined for $\xi<\la$.
Let $\cM_\la$ be the least model $\cM$ of $\zfc'$ of the form
$\rL_\ka\yt\ka<\omi$,
containing $\sis{\dU_\xi}{\xi<\la}$ and such that
$\la<\omi^\cM$ and all sets $\dU_\xi$, $\xi<\la$,
are countable in $\cM$.
Then $\dP_\la=\bigcup_{\xi<\la}\dU_\xi$ is countable in $\cM$, too.
Let $\sis{\vpi(j)}{j<\om}$ be the $\le_\rL$-least sequence of
\mut s $\vpi(j)\in\mt{\dP_\la}$, \dd\cle increasing and generic
over $\cM_\la$.
Define $\dU_\la=\dU$ as in Definition~\ref{dPhi}.
This completes the inductive step.

Let  $\dP=\bigcup_{\xi<\omi}\dU_\xi$.
\edf


\bpro
[in $\rL$]
\lam{uxip}
The sequence\/ $\sis{\dU_\xi}{\xi<\omi}$ belongs to\/
$\id{\hc}1$.\qed
\epro

\ble
[in $\rL$]
\lam{jden}
If a set\/ $D\in\cM_\xi\yt D\sq {\dP_\xi}$ is
pre-dense in\/ ${\dP_\xi}$ then it remains pre-dense in\/
$\dP$.
Therefore if\/ $\xi<\omi$ then\/
${\dU_\xi}$ is pre-dense in\/ $\dP$.

If a set\/ $D\in\cM_\xi\yt D\sq {\dP_\xi}\ti_{\Eo}{\dP_\xi}$
is pre-dense in\/ ${\dP_\xi}\ti_{\Eo}{\dP_\xi}$
then it is pre-dense in\/ $\dPd$.
\ele
\bpf
By induction on $\la\ge \xi$,
if $D$ is pre-dense in ${\dP_\la}$ then it
remains pre-dense in
${\dP_{\la+1}}=\dP_\la\cup\dU_\la$
by Lemma~\ref{uu3}.
Limit steps are obvious.
To prove the second claim note that
${\dU_\xi}$ is dense in ${\dP_{\xi+1}}$ by Lemma~\ref{uu2},
and $\dU_\xi\in\cM_{\xi+1}$.

To prove the last claim use Lemma~\ref{uu3d}.
\epf

\ble
[in $\rL$]
\lam{club}
If\/ $X\sq\HC=\rL_{\omi}$ then the set\/ $W_X$ of all
ordinals\/ $\xi<\omi$ such that\/
$\stk{\rL_\xi}{X\cap\rL_\xi}$ is an elementary submodel of\/
$\stk{\rL_{\omi}}{X}$ and\/ $X\cap\rL_\xi\in\cM_\xi$
is unbounded in\/ $\omi$.
More generally, if\/ $X_n\sq\HC$ for all\/ $n$
then the set\/ $W$ of all
ordinals\/ $\xi<\omi$, such that\/
$\stk{\rL_\xi}{\sis{X_n\cap\rL_\xi}{n<\om}}$
is an elementary submodel of\/
$\stk{\rL_{\omi}}{\sis{X_n}{n<\om}}$
and\/ $\sis{X_n\cap\rL_\xi}{n<\om}\in\cM_\xi$,
is unbounded in\/ $\omi$.
\ele
\bpf
Let $\xi_0<\omi$.
Let $M$ be a countable elementary submodel of $\rL_{\om_2}$
containing $\xi_0\yi\omi\yi X$,
and such that $M\cap\HC$ is transitive.
Let $\phi:M\onto\rL_\la$ be the Mostowski collapse, and
let $\xi=\phi(\omi)$.
Then $\xi_0<\xi<\la<\omi$ and $\phi(X)=X\cap\rL_\xi$ by the
choice of $M$.
It follows that $\stk{\rL_\xi}{X\cap\rL_\xi}$ is an elementary
submodel of $\stk{\rL_{\omi}}{X}$.
Moreover, $\xi$ is uncountable in $\rL_\la$, hence
$\rL_\la\sq\cM_\xi$.
We conclude that $X\cap\rL_\xi\in\cM_\xi$ since
$X\cap\rL_\xi\in\rL_\la$ by construction.

The second claim does not differ much: we start with a
model $M$ containing both the whole sequence
$\sis{X_n}{n<\om}$ and each particular $X_n$, and so on.
\epf

\bcor
[compare to \cite{jenmin}, Lemma 6]
\lam{ccc}
The forcing notions\/ $\dP$ and\/ $\dpd$ satisfy CCC in $\rL$.
\ecor
\bpf
Suppose that $A\sq\dP$ is a maximal antichain.
By Lemma~\ref{club}, there is an ordinal $\xi$ such that
$A'=A\cap{\dP_\xi}$ is a maximal antichain in ${\dP_\xi}$
and $A'\in\cM_\xi$.
But then $A'$ remains pre-dense, therefore,
still a maximal antichain, in the
whole set $\dP$ by Lemma~\ref{jden}.
It follows that $A=A'$ is countable.
\epf

\section{The model}
\las{mod}

We view the sets $\dP$ and $\dpd$
(Definition~\ref{uxi})
as forcing notions over $\rL$.

\ble
[compare to Lemma 7 in \cite{jenmin}]
\lam{mod1}
A real\/ $x\in\dn$ is\/ $\dP$-generic over\/ $\rL$ iff\/
$x\in Z=\bigcap_{\xi<\omi^\rL}\bigcup_{U\in\dU_\xi}[U]$.
\ele
\bpf
If $\xi<\omi^\rL$ then $\dU_\xi$ is pre-dense in $\dP$
by Lemma~\ref{jden}, therefore any real $x\in\dn$
$\dP$-generic over $\rL$ belongs to $\bigcup_{U\in\dU_\xi}[U]$.

To prove the converse, suppose that $x\in Z$ and prove that
$x$ is $\dP$-generic over $\rL$.
Consider a maximal antichain $A\sq\dP$ in $\rL$; we have to
prove that $x\in \bigcup_{T\in A}[T]$.
Note that $A\sq\dP_\xi$ for some
$\xi<\omi^\rL$ by Corollary~\ref{ccc}.
But then every tree $U\in\dU_\xi$ satisfies
$U\sqf\bigcup A$ by Lemma~\ref{uu3}, so that
$\bigcup_{U\in \dU_\xi}[U]\sq \bigcup_{T\in A}[T]$, and hence
$x\in \bigcup_{T\in A}[T]$, as required.
\epf

\bcor
[compare to Corollary 9 in \cite{jenmin}]
\lam{mod2}
In any generic extension of\/ $\rL$, the set of all reals
in\/ $\dn$ $\dP$-generic over\/ $\rL$ is\/ $\ip\HC1$ and\/
$\ip12$.
\ecor
\bpf
Use Lemma~\ref{mod1} and Proposition~\ref{uxip}.
\epf

\bdf
\lam{gg}
From now on, we assume that $G\sq\dPd$ is a set
\dd\dpdb generic over $\rL$,
so that the intersection
$X=\bigcap_{\ang{T,T'}\in G}[T]\ti[T']$ is
a singleton $X_G=\ans{\ang{x\lef[G],x\rig[G]}}$.
\edf

Compare the next lemma to Lemma 10 in \cite{jenmin}.
While Jensen's forcing notion in \cite{jenmin} guarantees
that there is a single generic real in the extension,
the forcing notion $\dP$ we use adds a whole \dd\Eo class
(a countable set) of generic reals!

\ble
[under the assumptions of Definition~\ref{gg}]
\lam{only}
If\/ $y\in\rL[G]\cap\dn$ then\/ $y$ is a\/ $\dP$-generic real
over\/ $\rL$ iff\/
$y\in\eko{x\lef[G]}\cup\eko{x\rig[G]}$.
\ele
Recall that $\eko x= \ens{\sg\app x}{\sg\in\bse}$.
\bpf
The reals $x\lef[G],\,x\rig[G]$ are separately \dd\dP generic
(see Remark~\ref{ref1}).
It follows that any real $y=\sg\app x\lef[G]\in\eko{x\lef[G]}$
or $y=\sg\app x\rig[G]\in\eko{x\rig[G]}$
is \dd\dP generic as well
since the forcing $\dP$ is
by definition invariant under the action of any $\sg\in\bse$.

To prove the converse, suppose towards the contrary that
there is a condition $\ang{T,T'}\in\dPd$ and
a \dd\dpdb real name $\rc=\sis{\kc ni}{n<\om,\,i=0,1}\in\rL$
such that $\ang{T,T'}$ \dd\dpdb forces that $\rc$ is
\dd\dP generic while $\dPd$ forces both formulas
$\rc\ne\sg\app\rpi\lef$ and $\rc\ne\sg\app\rpi\lef$ for all
$\sg\in\bse$.

Let $C_n=\kc n0\cup\kc n1$, this is a pre-dense set in $\dPd$.
It follows from Lemma~\ref{club} that there exists an ordinal
$\la<\omi$ such that
each set $C'_n=C_n\cap ({\dP_\la}\ti_{\Eo}{\dP_\la})$
is pre-dense in ${\dP_\la}\ti_{\Eo}{\dP_\la}$,
and the sequence $\sis{\xc ni}{n<\om,\,i=0,1}$ belongs to
$\cM_\la$, where $\xc ni=C'_n\cap \kc ni$ ---
then $C'_n$ is pre-dense in $\dPd$ too, by Lemma~\ref{jden}.
Therefore we can assume that in fact $C_n=C'_n$, that is,
$\rc\in\cM_\la$ and $\rc$ is a
\dd{({\dP_\la}\ti_{\Eo}{\dP_\la})}real name.

Further, as $\dpd$ forces that $\rc\ne\sg\app\rpi\lef$
and $\rc\ne\sg\app\rpi\rig$,
the set $D(\sg)$
of all conditions $\ang{S,S'}\in\dPd$ which directly force
$\rc\ne\sg\app\rpi\lef$ and $\rc\ne\sg\app\rpi\rig$,
is dense in $\dPd$ ---
for every $\sg\in\bse.$
Therefore, still by Lemma \ref{club}, we may
assume that the same ordinal $\la$ as above satisfies the
following: each set
$D'(\sg)=D(\sg)\cap{({\dP_\la}\ti_{\Eo}{\dP_\la})}$ is dense in
${\dP_\la}\ti_{\Eo}{\dP_\la}$.

Applying Theorem~\ref{K} with $\dP=\dP_\la$, $\dU=\dU_\la$,
and $\dP\cup\dU=\dP_{\la+1}$, we conclude that for each tree
$U\in\dU_\la$ the set $Q_U$ of all conditions
$\ang{V,V'}\in {\dP_{\la+1}}\ti_{\Eo}{\dP_{\la+1}}$
which directly force
$\rc\nin[U]$, is dense in ${\dP_{\la+1}}\ti_{\Eo}{\dP_{\la+1}}$.
As obviously $Q_U\in\cM_{\la+1}$, we further conclude that
$Q_U$ is pre-dense in the whole forcing $\dPd$
by Lemma~\ref{jden}.
This implies that $\dPd$ forces
$\rc\nin\bigcup_{U\in \dU_\la}[U]$,
hence, forces that $\rc$ is not \dd\dP generic,
by Lemma~\ref{mod1}.
But this contradicts to the choice of $\ang{T,T'}$.
\epf

\bcor
\lam{corf}
The set\/ $\eko{x\lef[G]}\cup\eko{x\rig[G]}$
is\/ $\ip12$ set in $\rL[G]$.
Therefore the 2-element set\/
$\ans{\eko{x\lef[G]},\eko{x\rig[G]}}$ is $\od$ in $\rL[G]$.\qed
\ecor

\bcor
\lam{cnr}
The\/ \dd\Eo classes\/
$\eko{x\lef[G]},\,\eko{x\rig[G]}$ are disjoint.
\ecor
\bpf
Corollary~\ref{lnr} implies\/ $x\lef[G]\not\Eo x\rig[G]$.
\epf

\vyk{
\ble
\lam{22}
The reals\/ $x\lef[G]$, $x\rig[G]$ are not\/ \dd\Eo equivalent,
therefore the\/ \dd\Eo classes\/
$\eko{x\lef[G]},\,\eko{x\rig[G]}$ are disjoint.
\ele
\bpf
Otherwise there is a string $\sg\in\bse$ and condition
$\ang{T,T'}\in \dpd$ which forces $\rpi\lef=\sg\app \rpi\rig$.
But it does not take much effort to define a stronger
condition $\ang{S,S'}\in \dpd$, $\ang{S,S'}\leq\ang{T,T'}$,
such that $[S]\cap(\sg\app[S'])=\pu$, easily leading to
contradiction.
\epf
}

\ble
[still under the assumptions of Definition~\ref{gg}]
\lam{sym}
Neither of the two\/ \dd\Eo classes\/
$\eko{x\lef[G]},\,\eko{x\rig[G]}$ is\/ $\od$ in\/ $\rL[G]$.
\ele
\bpf
Suppose towards the contrary that there is a
condition $\ang{T,T'}\in G$ and a formula $\vt(x)$ with
ordinal parameters such that $\ang{T,T'}$ \dd\dpdb forces
that $\vt(\eko{\rpi\lef})$ but $\neg\:\vt(\eko{\rpi\rig})$.
However both the formula and the forcing are invariant under
actions of strings in $\bse.$
In particular if $\sg\in\bse$ then
$\ang{\sg\app T,\sg\app T'}$ still \dd\dpdb forces
$\vt(\eko{\rpi\lef})$ and $\neg\:\vt(\eko{\rpi\rig})$.
We can take $\sg$ which satisfies $T'=\sg\app T$; thus
$\ang{T',T}$ still \dd\dpdb forces
$\vt(\eko{\rpi\lef})$ and $\neg\:\vt(\eko{\rpi\rig})$.%
\snos
{This is the argument which does not go through for the
full product $\dP\ti\dP$.}
However $\dpd$ is symmetric with respect to the
left-right exchange,
which implies that conversely $\ang{T',T}$ has to force
$\vt(\eko{\rpi\rig})$ and $\neg\:\vt(\eko{\rpi\lef})$.
The contradiction proves the lemma.
\epf
\qeD{Theorem~\ref{mt}}

\section{Conclusive remarks}
\las{konk}

(I)
One may ask whether other Borel \eqr s $\rE$ admit results
similar to Theorem~\ref{mt}.
Fortunately this question can be easily solved on the base
of the Glimm -- Effros dichotomy theorem \cite{hkl}.

\bcor
\lam{mtc}
The following is true in the model of Theorem~\ref{mt}.
Let\/ $\rE$ be a Borel \eqr\ on $\bn$ coded in\/ $\rL$.
Then  there  exists an\/ $\od$ pair of\/
\dd\rE equivalence classes\/ $\ans{\eke x,\eke y}$ such that
neither of the classes\/ $\eke x,\,\eke y$ is separately\/ $\od$,
iff \/ $\rE$ is not smooth.
\ecor
\bpf
Suppose first that $\rE$ is smooth.
By the Shoenfield absoluteness theorem,
the smoothness can be witnessed by a Borel
map $\vt:\bn\to\bn$ coded in $\rL$, hence, $\vt$ is $\od$ itself.
If $p=\ans{\eke x,\eke y}$ is $\od$ in the extension then so is
the 2-element set $R=\ens{\vt(z)}{z\in \eke x\cup\eke y}\sq\bn,$
whose both elements (reals), say $p_x$ and $p_y$,
are $\od$ by obvious reasons.
Then finally $\eke x=\vt^{(-1)}(p_x)$ and
$\eke y=\vt^{(-1)}(p_y)$
are $\od$ as required.

Now let $\rE$ be non-smooth.
Then by Shoenfield and the Glimm -- Effros dichotomy theorem
in \cite{hkl}, there is a continuous, coded by some
$r\in\bn\cap\rL$, hence, $\od$,
reduction $\vt:\dn\to\bn$ of $\Eo$ to $\rE$, so that we have
$a\Eo b$ iff $\vt(a)\rE\vt(b)$ for all $a,b\in\dn.$
Let, by Theorem~\ref{mt}, $\ans{\eko a,\eko b}$ be a $\ip12$ pair
of non-$\od$ \dd\Eo equivalence classes.
By the choice of $\vt$, one easily proves that
$\ans{\eke{\vt(a)},\eke{\vt(b)}}$
is a $\ip12(r)$ pair of non-$\od$ \dd\rE equivalence classes.
\epf

(II)
One may ask what happens with the Groszek -- Laver pairs
of sets of reals in better known models.
For some of them the answer tends to be in the negative.
Consider \eg\ the Solovay model of $\zfc$ in which all
projective sets of reals are Lebesgue measurable \cite{sol}.
Arguing in the Solovay model, let $\ans{X,Y}$ be an $\od$ set,
where $X,Y\sq\dn.$
Then the set of \rit{four} sets
$X\bez Y,\,Y\bez X,\,X\cap Y,\,\dn\bez(X\cup Y)$
is still $\od$, and hence we have an $\od$ \eqr $\rE$ on $\dn$
with four (or fewer if say $X\sq Y$) equivalence classes.
By a theorem of \cite{k:sol}\snos
{\label{f6}%
To replace the following brief argument,
one can also refer to a result by Stern implicit in
\cite{stl}: in the Solovay model, if an $\od$ \eqr\ $\rE$ has
at least one non-OD equivalence class then there is a pairwise
\dd\rE inequivalent perfect set.}
, either $\rE$ admits an $\od$
reduction $\vt:\dn\to2^{<\omi}$
to equality on $2^{<\omi}$ or $\Eo$ admits a
continuous reduction to $\rE$.
The ``or'' option fails since $\rE$ has finitely
many classes.

The ``either'' option leads to a finite
(not more than 4 elements) $\od$ set
$R=\ran\vt\sq2^{<\omi}$.
An easy argument shows that then every $r\in R$ is $\od$, and
hence so is the corresponding \dd\rE class $\vt^{-1}(r)$.
It follows that $X,Y$ themselves are $\od$.

\bvo
Is it true in the Solovay model that every \rit{countable}
$\od$ set $W\sq\pws{\bn}$ \rit{of sets of reals}
contains an $\od$ element $X\in W$ (a set of reals)?
\evo

An uncountable counterexample readily exists, for take the set
of all non-OD sets of reals.
As for sets $W\sq\bn$, any countable $\od$ set of reals in
the Solovay model consists of $\od$ elements, \eg\ by the
result mentioned in Footnote~\ref{f6}.

(III)
One may ask whether a forcing
similar to $\dpd$ with respect to the results in
Section~\ref{mod}, exists in ground models other than
$\rL$ or $\rL[x]$, $x\in\dn.$
Some coding forcing constructions with perfect trees
do exist
in such a general frameworks, see \cite{k-bag,kl:smz}.

\back
The authors thank Ali Enayat for the interest
in the problem and helpful remarks.
The authors thank the anonymous referee for many
important suggestions that helped to improve the text.
\eack

\bibliographystyle{amsplain}

\begin{thebibliography}{10}

\bibitem{k-bag}
J.~Bagaria and V.~Kanovei,
\newblock
\textit{On coding uncountable sets by reals},
\newblock
Mathematical Logic Quarterly,
56(4): 409--424, 2010.

\bibitem{ena}
Ali {Enayat},
\newblock
\textit{On the Leibniz-Mycielski axiom in set theory},
\newblock {Fundam. Math.}, 181(3): 215--231, 2004.

\vyk{
\bibitem{sdf}
Sy D.~Friedman.
\newblock {The $\Pi^1_2$-singleton conjecture.}
\newblock {\em {J. Amer. Math. Soc.}}, 3(4): 771--791, 1990.
}

\bibitem{gl}
M.~{Groszek} and R.~{Laver},
\newblock \textit{Finite groups of OD-conjugates},
\newblock {Period. Math. Hung.}, 18: 87--97, 1987.

\bibitem{hkl}
L.~A. Harrington, A.~S. Kechris, and A.~Louveau,
\newblock
\textit{A {G}limm-{E}ffros dichotomy for
{B}orel equivalence relations},
\newblock {J. Amer. Math. Soc.}, 3(4): 903--928, 1990.

\bibitem{hms}
L.~A. Harrington, D. Marker, and S.~Shelah,
\newblock \textit{{B}orel orderings},
\newblock {Trans.\ Amer.\ Math.\ Soc.}, 310(1): 293--302, 1988.

\bibitem{jechmill}
Thomas {Jech},
\newblock \textit{Set theory},
\newblock Berlin: Springer,
the third millennium revised and expanded edition, 2003.

\bibitem{jenmin}
Ronald {Jensen},
\newblock
\textit{Definable sets of minimal degree},
\newblock {Math. Logic Found. Set Theory, Proc. Int. Colloqu.,
Jerusalem 1968,
  pp.\ 122-128}, 1970.


\bibitem{k:sol}
Vladimir Kanovei,
\newblock
\textit{An {U}lm-type classification theorem for
equivalence relations in {S}olovay model},
\newblock {J. Symbolic Logic}, 62(4): 1333--1351, 1997.

\bibitem{kanB}
Vladimir Kanovei,
\newblock \textit{Borel equivalence relations.
Structure and classification},
\newblock Providence, RI:
American Mathematical Society (AMS), 2008.

\bibitem{kl:smz}
V.~{Kanovei} and V.~{Lyubetsky},
\newblock
\textit{An effective minimal encoding of uncountable sets},
\newblock
Siberian Mathematical Journal, 52(5): 854--863, 2011.


\bibitem{kl:cds}
V.~{Kanovei} and V.~{Lyubetsky},
\newblock \textit{A countable definable set of reals
containing no definable elements},
\newblock {ArXiv e-prints, 1408.3901}, August 2014.

\bibitem{kl:ceo}
V.~{Kanovei} and V.~{Lyubetsky},
\newblock \textit{A definable \dd{E_0}class containing no
definable elements},
\newblock
Archive for Mathematical Logic, 2015, 54, 5, pp. 711--723.

\bibitem{ksz}
Vladimir Kanovei, Martin Sabok, and Jind\v{r}ich Zapletal,
\newblock \textit{Canonical Ramsey theory on Polish space},
\newblock Cambridge: Cambridge University Press, 2013.

\bibitem{sol}
R.M. Solovay,
\newblock \textit{A model of set-theory in which every set
of reals is Lebesgue measurable},
\newblock {Ann. Math. (2)}, 92: 1--56, 1970.

\bibitem{stl}
J. Stern,
\newblock \textit{On Lusin's restricted continuum problem},
\newblock {Ann. Math. (2)}, 120: 7--37, 1984.


\end{thebibliography}

\end{document}